\documentclass[10pt, a4paper]{amsart}
\usepackage{amsmath}
\usepackage{amssymb}
\usepackage{latexsym}
\usepackage{a4}
\usepackage[all]{xy}

\usepackage{hyperref}

\theoremstyle{plain}
\newtheorem*{thm*}{Theorem}
\newtheorem*{rem*}{Remark}

\newtheorem{thm}{Theorem}[section]
\newtheorem{cor}[thm]{Corollary}
\newtheorem{defi}[thm]{Definition}
\newtheorem{prop}[thm]{Proposition}
\newtheorem{lm}[thm]{Lemma}

\newtheorem{claim*}{Claim}
\newtheorem{exple}[thm]{Example}
\newtheorem{rem}[thm]{Remark}
\newtheorem{nota}[thm]{Notation}
\numberwithin{equation}{thm}

\newcommand{\A}{{\mathcal{A}}}
\newcommand{\C}{{\mathcal{C}}}
\newcommand{\D}{{\mathcal{D}}}
\newcommand{\F}{{\mathcal{F}}}
\newcommand{\E}{{\mathcal{E}}}
\newcommand{\FF}{{\mathbb{F}_2}}

\newcommand{\Eq}{{\ensuremath{\mathcal{E}_{q}}}}
\newcommand{\Eqd}{{\mathcal{E}_{q}^{\mathrm{deg}}}}
\newcommand{\Tq}{{\ensuremath{\mathcal{T}_{q}}}}
\newcommand{\Sq}{{\ensuremath{\mathcal{S}_{q}}}}
\newcommand{\Fq}{{\ensuremath{\mathcal{F}_{quad}}}}
\newcommand{\Gq}{{\ensuremath{\mathcal{F}_{iso}}}}

\begin{document}

\title{The functor category $\mathcal{F}_{quad}$}

\author{Christine Vespa}

\address{Laboratoire Analyse, G\'{e}om\'{e}trie et Applications, UMR 7539\\
  Institut Galil\'{e}e, Universit\'{e} Paris 13, 93430 Villetaneuse, France}
\email{vespa@math.univ-paris13.fr}

\date{\today}

\begin{abstract}
In this paper, we define the functor category $\Fq$ associated to
$\FF$-vector spaces equipped with a quadratic form. We show the
existence of a fully-faithful, exact functor $\iota: \F \rightarrow
\Fq$, which preserves simple objects, where $\F$ is the category of
functors from the category of finite dimensional $\FF$-vector spaces to the category of all
$\FF$-vector spaces. We define the subcategory $\Gq$ of $\Fq$, which is
equivalent to the product of the categories of modules over the
orthogonal groups; the inclusion is a fully-faithful functor $\kappa: \Gq \rightarrow
\Fq$ which preserves simple objects.\\
$\ $
\end{abstract}

\maketitle


\section*{Introduction}

In recent years, one of the functor categories which has been
particularly studied is the category $\F(p)$ of functors from the category
$\E^f$ of finite dimensional $\mathbb{F}_p$-vector spaces to the category $\E$ of
all $\mathbb{F}_p$-vector spaces, where $\mathbb{F}_p$ is the prime
field with $p$ elements. This category is connected
to several areas of algebra and some examples of these can be
found in \cite{PS16}. The category $\F(p)$ is closely related to the
general linear groups. An important application of $\F(p)$ is given in
\cite{FFSS}, where the four authors proved that this category is very useful
for the study of the stable cohomology of the general linear
groups with suitable coefficients. They showed that the calculation  of certain extension
groups in the category $\F$  determines some stable cohomology groups
of general linear groups. One of the motivations of the work presented
here is to construct and study a category $\Fq$ which could 
play a similar role for the stable cohomology of the  orthogonal
groups.

In this paper, we restrict to the prime $p=2$; the techniques can be
applied in the odd prime case, but the case $p=2$ presents features
which make it particularly interesting. Henceforth, we will suppose
that $p=2$ and we will denote the category $\F(2)$ by
$\F$.

After some recollections on the theory of quadratic forms over $\FF$, we give
the definition of the category $\Fq$. In order to have a good
understanding of the category $\Fq$, we seek to classify the simple
objects of this category. The first important result of
this paper is the following theorem.
\begin{thm*}
There is a functor $$\iota: \F \rightarrow \Fq$$ which satisfies the
following properties:
\begin{enumerate}
\item
$\iota$ is exact;
\item
$\iota$ preserves tensor products;
\item
$\iota$ is fully-faithful;
\item
if $S$ is a simple object of $\F$, $\iota(S)$ is a simple object of $\Fq$.
\end{enumerate} 
\end{thm*}

\begin{rem*}
Using a study of the projective functors of the category $\Fq$, we
will show in \cite{Vespa} that $\iota(F)$ is a thick
subcategory of $\Fq$.
\end{rem*}

To study a particular family of functors of $\Fq$, the
isotropic functors, we define the subcategory $\Gq$ of $\Fq$ which is
related to $\Fq$ by the following theorem:

\begin{thm*}
There is a functor $$\kappa: \Gq \rightarrow \Fq$$ which satisfies the
following properties:
\begin{enumerate}
\item
$\kappa$ is exact;
\item
$\kappa$ preserves tensor products;
\item
$\kappa$ is fully-faithful;
\item
if $S$ is a simple object of $\Gq$, $\kappa(S)$ is a simple object of $\Fq$.
\end{enumerate} 
\end{thm*}

We obtain the classification of the simple objects of the category
$\Gq$ from the following theorem.
\begin{thm*}
There is a natural equivalence of categories 
$$\Gq \simeq \prod_{V \in \mathcal{S}} \FF[O(V)]-mod$$
where $\mathcal{S}$ is a set of representatives of isometry classes of
quadratic spaces (possibly degenerate).
\end{thm*}

Apart from the previous theorem, the results of this paper are contained in
the Ph.D. thesis of the author \cite{Vespa-these}, although several results
are presented here from a more conceptual point of view.

The author wants to thank her PhD supervisor, Lionel Schwartz, as well
as Geoffrey Powell and Aur\'{e}lien Djament for their useful comments
and suggestions on a previous version of this paper. 
\section{Quadratic spaces over $\FF$}

We recall the definition and the classification of quadratic forms over the field $\FF$. We refer the reader to
\cite{Pfister} for details.
\subsection{Definitions}
Let $V$ be a $\FF$-vector space of finite dimension. A quadratic form over $V$ is a function $q: V \rightarrow \FF$ such
that $q(x+y)+q(x)+q(y)=B(x,y)$ is a bilinear form.
As a direct consequence of the definition, we have that the bilinear form
associated to a quadratic form is alternating.

The radical of a quadratic space $(V,q_V)$ is the subspace of $V$ given by
$$\mathrm{Rad}(V,q_V)=\{ v \in V | \forall w \in V\  B(v,w)=0 \}$$
where $B(-,-)$ is the bilinear form associated to the quadratic form.
A quadratic space $(V,q_V)$ is non-degenerate if $\mathrm{Rad}(V,q_V)=0$.

\subsection{Non-degenerate quadratic forms}
\subsubsection{Classification}
In this paragraph, we recall the classification of non-degenerate
quadratic forms. The classification of non-singular
alternating forms implies that a non-degenerate quadratic space over
$\FF$ has even dimension and has a symplectic basis.

The space $H_0$ is the non-degenerate quadratic space of dimension two
with symplectic basis $\{ a_0, b_0 \}$, and quadratic
form determined by:
$$ \begin{array}{cccc}
q_{0}: & H_0 & \rightarrow & \mathbb{F}_{2}\\
       & a_0 & \longmapsto & 0\\
       & b_0 & \longmapsto & 0.
\end{array}$$
The space $H_{1}$ is the non-degenerate quadratic space of dimension two
with symplectic basis $\{ a_1, b_1 \}$, and quadratic
form determined by:
$$ \begin{array}{cccc}
q_{1}: & H_1 & \rightarrow & \mathbb{F}_{2}\\
       & a_1 & \longmapsto & 1\\
       & b_1 & \longmapsto & 1.
\end{array}$$
The spaces $H_0$ and $H_1$ are not isometric, whereas the spaces  $H_0
\bot H_0$ and $H_1 \bot H_1$ are isometric. The non
degenerate quadratic spaces of dimension $2m$, for $m \ge 1$, are
classified by the following result.
\begin{prop} \cite{Pfister} \label{1.1}
Let $m \geq 1$ be an integer.
\begin{enumerate}
\item
The quadratic spaces $H_0^{\bot m}$ and $H_1 \bot H_0^{\bot (m-1)}$ are not isometric.
\item
A quadratic space of dimension $2m$ is isometric to either $H_0^{\bot m}$ or $H_1 \bot H_0^{\bot (m-1)}$.
\end{enumerate}
\end{prop}
The two spaces $H_0^{\bot m}$ and $H_1 \bot H_0^{\bot (m-1)} $ are distinguished by the Arf
invariant introduced in \cite{Arf}. Observe that $\mathrm{Arf}(H_0^{\bot m})=0$ and $\mathrm{Arf}(H_1 \bot H_0^{\bot (m-1)})=1$.

\subsubsection{The category $\Eq$}
\begin{defi}
Let $\Eq$ be the category having as objects finite dimensional $\FF$-vector spaces equipped with a non
degenerate quadratic form and with morphisms linear maps which
preserve the quadratic forms.
\end{defi}
Observe that a  linear map which preserves the quadratic
forms preserves the underlying bilinear form, but the converse is, in
general, false. The following proposition summarizes straightforward but important results about the
category $\Eq$.
\begin{prop} \label{1.3}
\begin{enumerate}
\item
The morphisms of $\Eq$ are injective linear maps. Consequently they
are monomorphisms.
\item
The category $\Eq$ does not admit push-outs or pullbacks.
\end{enumerate}
\end{prop}
\begin{exple}
The diagram $V \leftarrow \{0 \} \rightarrow V$, where $V \not \simeq
\{ 0\}$ does not admit a push-out in $\Eq$.
\end{exple}
To resolve this difficulty, we define the notion of a pseudo push-out
in $\Eq$.

\subsubsection{Pseudo push-out}
To define the pseudo push-out in $\Eq$, we need the following remark,
which uses the non-degeneracy of the quadratic form in an essential way.
\begin{rem}
For $f$ an element of $\mathrm{Hom}_{\Eq}(V,W)$, we have $W \simeq
f(V) \bot V'$, where $V'$ is the orthogonal space to $W$ in the space $V$. As the spaces $V$ and $f(V)$ are isometric, the spaces
$W$ and $V \bot V'$ are also. We will write 
$$f: V \rightarrow W \simeq V \bot V'.$$
\end{rem}
\begin{defi} \label{1.6}
For $f: V \rightarrow W \simeq V \bot V'$ and $g: V \rightarrow X
\simeq V \bot V''$  morphisms in $\Eq$, the pseudo push-out of $f$ and
$g$, denoted by $W \underset{V}{\bot} X$,  is the object $V \bot V' \bot V''$ of $\Eq$.
\end{defi}
We give, in the following proposition, the principal properties of the
pseudo push-out.
\begin{prop} \label{1.7}
Let $f:V \rightarrow W \simeq V \bot V'$ and $g: V \rightarrow X
\simeq V \bot V'' $ be morphisms in $\Eq$, the pseudo push-out of $f$
and $g$ satisfies the following properties.
\begin{enumerate}
\item
There exists a commutative diagram of the form
 $$\xymatrix{
  V \ar[d]_f  \ar[r]^g  & X \ar[d]\\
  W \ar[r]   & X \underset{V}{\bot} W}$$
\item
$W \underset{V}{\bot} X \simeq X \underset{V}{\bot} W$;
\item
if $V \simeq V'$ then  $ W \underset{V}{\bot} X \simeq W
\underset{V'}{\bot} X $;
\item{Associativity:}
$(X \underset{V}{\bot} W) \underset{Z}{\bot} Y\simeq X
\underset{V}{\bot} (W \underset{Z}{\bot} Y)$;
\item{Unit: }
 $V \underset{V}{\bot} W \simeq W$.

\end{enumerate}
\end{prop}
\begin{rem}
By the first point of the previous proposition, the pseudo push-out
occurs in a commutative diagram
 $$\xymatrix{
  V \ar[d]  \ar[r]  & V \bot V" \ar[d]\\
  V \bot V' \ar[r]   & V \bot V' \bot V''}$$
which is equivalent to the orthogonal sum of the diagram 
 $$\xymatrix{
  \{0 \} \ar[d]  \ar[r]  & V'' \ar[d]\\
  V' \ar[r]   & V' \bot V''}$$
with $V$. Hence, the pseudo push-out can be considered as a
generalization of the orthogonal sum.
\end{rem}

\subsection{Degenerate quadratic forms}
The previous section implies that, in particular, all quadratic
spaces of odd dimension are degenerate. We begin by considering the
quadratic spaces of dimension one.
\begin{nota}
For $\alpha \in \{ 0,1 \}$, let $(x, \alpha)$ be  the quadratic space of dimension one generated by $x$ such that $q(x)=\alpha$.
\end{nota}
\subsubsection{Classification}
We have the following classification:
\begin{thm} \cite{Pfister} \label{1.9}
\begin{enumerate}
\item
Every quadratic space over $\FF$ has an orthogonal decomposition 
$$V \simeq H \bot \mathrm{Rad}(V)$$
such that $H$ is non-degenerate and $\mathrm{Rad}(V)$ is isometric to
either $(x,0)^{\bot r}$ or $(x,1)^{\bot r}$, where $r$ is the dimension of $\mathrm{Rad}(V)$.
\item
Let $V \simeq H \bot \mathrm{Rad}(V)$ and $V \simeq H' \bot
\mathrm{Rad}(V)$ be two decompositions of $V$, if $\mathrm{Rad}(V)
\simeq (x,0)^{\bot r}$ for $r \geq 0$, then $H \simeq H'$.
\item
Let $H$ and $H'$ be two non-degenerate quadratic forms such that
$\mathrm{dim}(H) =\mathrm{dim}(H')$ then, for all $r>0$,
$$H \bot (x,1)^{\bot r} \simeq H' \bot (x,1)^{\bot r}.$$
\end{enumerate}
\end{thm}

\begin{rem}
The third point of the previous theorem, is implied by the isometry:
$$ H_0 \bot (x,1) \simeq H_1 \bot (x,1).$$
This exhibits one of the particularities of quadratic forms over $\FF$:
the ``non-degenerate part'' of a quadratic form is not unique in general, not
even up to isometry.
\end{rem}
\subsubsection{The category $\Eqd$}
\begin{defi}
Let $\Eqd$ be the category having as objects finite dimensional $\FF$-vector spaces
equipped with a (possibly degenerate) quadratic
form and with morphisms, injective linear maps which
preserve the quadratic forms.
\end{defi}
\begin{rem}
The hypothesis that the morphisms are linear injective maps is
essential for later considerations.
\end{rem}
The following proposition underlines one of the important difference
between $\Eq$ and $\Eqd$.
\begin{prop} \label{1.12}
The category $\Eqd$ admits pullbacks.
\end{prop}
Nevertheless, $\Eqd$ does not contain all push-outs.

\section{Definition of the category $\Fq$}
We have emphasized the fact that all the morphisms of $\Eq$ are
monomorphisms. The constructions which link the category $\F$ to
the stable homology of the general linear groups use, in an essential
way, the existence of retractions in the category $\E^f$. Therefore, to
consider analogous constructions in the quadratic case, we would like
to add formally retractions to the category $\Eq$. For this, we
define the category $\Tq$ inspired by the construction of the category of
co-spans introduced by B\'{e}nabou in \cite{Benabou}.

\subsection{The category $\Tq$}

\subsubsection{The categories of spans and co-spans}
We will begin by recalling the construction of B\'{e}nabou.
\begin{rem}
Our principal interest is in the category $\Fq$. The construction of
this category uses a generalization of the category of co-spans, which
we have chosen to present rather than spans.
\end{rem}

\begin{defi} \cite{Benabou}
Let $\D$ be a category equipped with push-outs, the category
$\mathrm{coSpan}(\D)$ is defined in the following way:
\begin{enumerate}
\item
the objects of $\mathrm{coSpan}(\D)$ are those of $\D$;
\item
for $A$ and $B$ two objects of $\mathrm{coSpan}(\D)$,
$\mathrm{Hom}_{\mathrm{coSpan}(\D)} (A,B)$ is the set of equivalence
classes of diagrams in $\D$ of the form $A \xrightarrow{f} D \xleftarrow{g} B$, for the equivalence relation which
identifies the two diagrams $A \xrightarrow{f} D \xleftarrow{g} B$ and
$A \xrightarrow{u} D' \xleftarrow{v} B$ if there exists an isomorphism
$\alpha: D \rightarrow D'$ such that the following diagram is
commutative
$$\xymatrix{
  &  B \ar[d]_{g} \ar[rdd]^{v}\\
A \ar[r]^{f}  \ar[drr]_{u} & D \ar[dr]^(.3){\alpha }\\
  &   &D'.\\
}$$
The morphism of $\mathrm{Hom}_{\mathrm{coSpan}(\D)} (A,B)$ represented
by the diagram $A \xrightarrow{f} D \xleftarrow{g} B$ will be denoted
by $[A \xrightarrow{f} D \xleftarrow{g} B]$;
\item the composition is given by:\\
for two
  morphisms $T_1 =\lbrack A \rightarrow D  \leftarrow B
 \rbrack$ and  $T_2 =\lbrack B \rightarrow D' \leftarrow C \rbrack$,
$$T_2 \circ T_1\ =\lbrack A \rightarrow D \underset{B}{\amalg}
 D' \leftarrow C \rbrack.$$
\end{enumerate}
\end{defi}

By duality, B\'{e}nabou gives the following definition.

\begin{defi}
Let $\D$ be a category equipped with pullbacks. The category
$\mathrm{Sp}(\D)$ is defined by: $\mathrm{Sp}(\D) \simeq \mathrm{coSpan}(\D^{op})^{op}$.
\end{defi}
\begin{exple} \label{2.3}
The category $\mathrm{Sp}(\Eqd)$ is defined, since $\Eqd$ admits
pullbacks by proposition \ref{1.12}.
\end{exple}
\subsubsection{The category $\Tq$}
By proposition \ref{1.3} neither the category of Spans nor co-Spans
are defined for $\Eq$. However, we observe that the universality of the
push-out plays no role in the definition of the category
$\mathrm{coSpan}(\D)$. So, by definition \ref{1.6} of the
pseudo push-out of $\Eq$ we can give the following definition.
\begin{defi}
The category $\hat{\mathcal{T}_{q}}$ is defined in the following way:
\begin{enumerate}
\item
the objects of $\hat{\mathcal{T}_{q}}$ are those of $\Eq$;
\item for $V$ and $W$ two objects of $\hat{\mathcal{T}_{q}}$,
$\mathrm{Hom}_{\hat{\mathcal{T}_{q}}}(V,W)$ is defined in the same way
as for the category $\mathrm{coSpan}(\D)$ and we will use the same
notation for the morphisms;
\item
the composition is given by:\\
for $T_1 =\lbrack V \rightarrow X_1 \leftarrow W
   \rbrack$ and $T_2 =\lbrack W \rightarrow X_2 \leftarrow Y \rbrack$,
$$T_2 \circ T_1\ =\lbrack V \rightarrow X_1
   \underset{W}{\bot} X_2 \leftarrow Y \rbrack$$ 
where $ X_1
   \underset{W}{\bot} X_2$ is the pseudo push-out.
\end{enumerate}
\end{defi}
The elementary properties of the pseudo push-out given in \ref{1.7}
show that the composition is well-defined and associative. Thus, the
above defines a category.

To define the category $\Tq$, we consider the following relation on
the morphisms of $\hat{\Tq}$.

\begin{defi}
For $V$ and $W$ two objects of $\hat{\Tq}$, the relation $\mathcal{R}$
on $\mathrm{Hom}_{\hat{\mathcal{T}_{q}}}(V,W)$  for $V$ and $W$
 objects of $\hat{\mathcal{T}_{q}}$, is defined by:\\
 for $T_1 =\lbrack V \rightarrow X_1 \leftarrow W
   \rbrack$ and $T_2 =\lbrack V \rightarrow X_2 \leftarrow W \rbrack$,
   two elements of $\mathrm{Hom}_{\hat{\mathcal{T}_{q}}}(V,W)$, $T_1
   \mathcal{R} T_2$ if there exists a morphism 
$\alpha$  of $\Eq$\ such that the following diagram is commutative
$$\xymatrix{
  & W \ar[d] \ar[rdd]\\
V \ar[r]  \ar[drr] & X_{1}\ar[dr]^(.3){\alpha }\\
  &   &X_{2}.\\
}$$
We will denote by $\sim$ the equivalence relation on
$\mathrm{Hom}_{\hat{\mathcal{T}_{q}}}(V,W)$ generated by the relation
$\mathcal{R}$.
\end{defi}

\begin{lm}
The composition in $\hat{\mathcal{T}_{q}}$ induces an application:
$$ \circ:\ \mathrm{Hom}_{\hat{\mathcal{T}_{q}}}(V,W)/ \sim \  \times \ 
\mathrm{Hom}_{\hat{\mathcal{T}_{q}}}(U,V) / \sim \ \rightarrow\ 
\mathrm{Hom}_{\hat{\mathcal{T}_{q}}}(U,W)/ \sim.$$
\end{lm}
\begin{proof}
By the properties of the pseudo push-out given in Proposition
\ref{1.7} we verify that:
\begin{enumerate}
\item
if $T_1 \mathcal{R} T_2$, then $(T_3 \circ T_1)
\mathcal{R} (T_3 \circ T_2)$;
\item
if  $T_3 \mathcal{R} T_4$, then $(T_1 \circ T_3)
\mathcal{R} (T_1 \circ T_4)$.
\end{enumerate}
\end{proof}
Thanks to the previous lemma, we can give the following definition.

\begin{defi} 
Let $\Tq$ be the category having as objects the objects of $\Eq$ and
with morphisms $  \mathrm{Hom}_{{\mathcal{T}_{q}}}(V,W)=\mathrm{Hom}_{\hat{\mathcal{T}_{q}}}(V,W)/\sim$.
\end{defi}

For convenience we will use the same notation for the morphisms of
$\Tq$ as for those of $\hat{\mathcal{T}_{q}}$.

\subsubsection{Properties of the category  $\Tq$} \label{2.1.3}
We have the following important property.
\begin{prop}
For $f: V \rightarrow W$ a morphism in the category $\Eq$, we have
the following relation in $\Tq$
$$[W \xrightarrow{\mathrm{Id}} W \xleftarrow{f} V ] \circ [V
\xrightarrow{f} W \xleftarrow{\mathrm{Id}} W] =\mathrm{Id}_V.$$
In particular, $[W \xrightarrow{\mathrm{Id}} W \xleftarrow{f} V ]$ is a
retraction of $[V \xrightarrow{f} W \xleftarrow{\mathrm{Id}} W]$.
\end{prop} 
\begin{proof}
It is a direct consequence of the definition of $\Tq$.
\end{proof}
To close this paragraph we give two useful constructions in the
category $\Tq$.

\begin{defi}[Transposition] \label{2.9}
The transposition functor, $tr: \Tq^{op}
\rightarrow \Tq$, is defined on objects by $tr(V)=V$ and on morphisms by:
$$\mathrm{tr}(f)=\mathrm{tr} (\lbrack  V \rightarrow X \leftarrow W
\rbrack)= \lbrack  W \rightarrow X \leftarrow V \rbrack,$$
for $V$ and $W$ objects of $\Tq^{op}$ and $f$
an element of $\mathrm{Hom}_{\Tq^{op}}(W,V)$.
\end{defi}
Observe that the transposition functor is involutive.
\begin{prop}[Orthogonal sum] \label{2.10}
There exists a bifunctor $\bot :\ \Tq \times \Tq \rightarrow \Tq$,
called the orthogonal sum, defined on objects by:
$$\bot (V,W)=V \bot W$$
and on morphisms by:
$$\bot (\lbrack V
\xrightarrow{f} X \xleftarrow{g} W \rbrack, \lbrack V'
\xrightarrow{f'} X' \xleftarrow{g'} W' \rbrack)=\lbrack V \bot V'
\xrightarrow{f \bot f'} X \bot X' \xleftarrow{g \bot g'} W \bot W'
\rbrack.$$
This bifunctor gives $\Tq$ the structure of a symmetric monoidal
category, with unit $\{0 \}$.
\end{prop}
\begin{proof}
It is straightforward to verify that $\bot$ is a well-defined
bifunctor, and that it is associative,
symmetric and that the object $\{ 0 \}$ of $\Tq$ is a unit for $\bot$.
\end{proof}

\subsection{Definition and properties of the category $\Fq$}
As the category $\Tq$ is essentially small, we can give the following
definition.
\begin{defi}
The category $\Fq$ is the category of functors from $\Tq$ to $\E$.
\end{defi}

\begin{rem}
By analogy with the classical definition of Mackey functors, given in
\cite{thev} and due, originally, to Dress \cite{dress}, and with the
work of Lindner in \cite{lindner}, we can view the category $\Fq$ as
the category of generalized Mackey functors over $\Eq$. Note that,
in our definition, we consider all the functors from $\Tq$ to $\E$ and
not only the additive functors, unlike \cite{thev}.
\end{rem}

By classical results about functor categories and by the Yoneda lemma,
we obtain the following theorem.

\begin{thm} \label{2.13}
\begin{enumerate}
\item
The category $\Fq$ is abelian.
\item
The tensor product of vector spaces induces a structure of symmetric
monoidal category on $\Fq$.
\item
For any object $V$ of $\Tq$, the functor $P_V=\FF \lbrack
\mathrm{Hom}_{\Tq}(V,-) \rbrack $ is a projective object and there is
a natural isomorphism:
$$\mathrm{Hom}_{\Fq}(P_V, F) \simeq F(V)$$
for all objects $F$ of $\Fq$.

The set of functors $\{ P_V | V \in \mathcal{S} \}$ is a set of
projective generators of $\Fq$, where $\mathcal{S}$ is a set of
representatives of isometry classes of non-degenerate quadratic spaces. In particular, the category $\Fq$ has enough projective objects.
\end{enumerate}
\end{thm}
The transposition functor of $\Tq$ allows us to give the following definition.
\begin{defi} \label{2.14}
The duality functor of $\Fq$ is the functor $D: \Fq^{op} \rightarrow
\Fq$ given by:
$$DF=-^* \circ F \circ tr^{op}$$
for $F$ an object of $\Fq$, $-^*$ the duality functor from $\E^{op}$ to
$\E$ and $tr$ the transposition functor of $\Tq$ defined in \ref{2.9}.
\end{defi}

The following proposition summarizes the basic properties of the duality
functor $D$.

\begin{prop}
\begin{enumerate}
\item
The functor $D$ is exact.
\item
The functor $D$ is right adjoint to the functor $D^{op}$, i.e. we have
a natural isomorphism:
$$\mathrm{Hom}_{\Fq}(F, DG) \simeq \mathrm{Hom}_{\Fq^{op}}(D^{op}F, G) \simeq \mathrm{Hom}_{\Fq}(G, DF).$$
\item
For $F$ an object of $\Fq$ with values in finite dimensional vector
spaces, the unit of the adjunction between $\Fq$ and $\Fq^{op}$, $F
\rightarrow DD^{op}F$, is an isomorphism.

\end{enumerate}
\end{prop}
A straightforward consequence of the second point of the last
proposition and theorem \ref{2.13} is:
\begin{cor}
The category $\Fq$ has enough injective objects.
\end{cor}

\section{Connection between $\F$ and $\Fq$}
Recall that $\F$ is the category of functors from $\E^f$ to $\E$, where $\E^f$ is the full
subcategory of $\E$ having as objects the finite dimensional spaces. The main result of this section is the following theorem.

\begin{thm} \label{3.1}
There is a functor $$\iota: \F \rightarrow \Fq$$ which satisfies the
following properties:
\begin{enumerate}
\item
$\iota$ is exact;
\item
$\iota$ preserves tensor products;
\item
$\iota$ is fully-faithful;
\item
if $S$ is a simple object of $\F$, $\iota(S)$ is a simple object of $\Fq$.
\end{enumerate} 
\end{thm}
To define the functor $\iota$ of the last theorem we need to define
the forgetful functor $\epsilon: \Tq \rightarrow \E^f$ which can be
viewed as an object of the category $\Fq$, by composing with $\E^f
\hookrightarrow \E$.
\subsection{The forgetful functor $\epsilon$ of $\Fq$}
\subsubsection{Definition}
\begin{nota}
We denote by $\mathcal{O}: \Eq \rightarrow \E^f$ the functor which
forgets the quadratic form.
\end{nota}
\begin{prop}
There exists a functor $\epsilon: \Tq \rightarrow \E^f$ defined by
$\epsilon(V)=\mathcal{O}(V)$ and 
$$\epsilon([V \xrightarrow{f} W \bot W' \xleftarrow{g} W ])=p_g \circ
\mathcal{O}(f)$$
where $p_g$ is the orthogonal projection from $W \bot W'$ to $W$.
\end{prop}
The proof of this proposition relies on the following
straightforward property of the pseudo push-out of $\Eq$.
\begin{lm}
For a pseudo push-out diagram in $\Eq$
 $$\xymatrix{
  V \ar[d]_f  \ar[r]^g  & V \bot V'' \ar[d]^{l}\\
  V \bot V' \ar[r]_-{k}    & V \bot V' \bot V''}$$
we have the following relation in $\E^f$:
$$\mathcal{O}(g) \circ p_f = p_l \circ \mathcal{O}(k)$$
where $p_f$ and $p_l$ are the orthogonal projections associated to,
respectively, $f$ and $l$.
\end{lm}
\begin{proof}[Proof of the proposition]
It is straightforward to check that the functor $\epsilon$ is well defined on
the classes of morphisms of $\Tq$. To show that $\epsilon$ is a
functor, we verify that 
$$\epsilon([V \xrightarrow{\mathrm{Id}} V \xleftarrow{\mathrm{Id}}
V])=\mathrm{Id}$$
and the relation for the composition is a direct consequence of the
above lemma.
\end{proof}

\subsubsection{The fullness of  $\epsilon$}
The aim of this section is to prove the following proposition:
\begin{prop}
The functor $\epsilon: \Tq \rightarrow \E^f$ is full.
\end{prop}
\begin{proof}
Let $(V,q_V)$ and $(W,q_W)$ be two objects of $\Tq$ and $f \in
\mathrm{Hom}_{\E^f}(\epsilon(V,q_V), \epsilon (W,q_W))$ be a linear
map from $V$ to $W$. We prove, by induction on the dimension of
$V$, that there is a morphism in $\Tq$: $T= [ V \xrightarrow{\varphi} X=W \bot
Y \xleftarrow{i} W]$ such that $\epsilon(T)=f$. The proof is based on
the idea that, for a sufficiently large space $X$, we can obtain all
the linear maps.

As the quadratic space $V$ is non-degenerate, we know that it has even
dimension.

To start the induction, let  $(V,q_V)$ be a non-degenerate quadratic space of dimension two, with
 symplectic basis $\{ a,b \}$ and  $f: V \rightarrow W$ be
a linear map. We verify that the following linear map:

$$ \begin{array}{cccl}
g_1: & V& \rightarrow & W \bot H_1 \bot H_0 \simeq W \bot
\mathrm{Vect(a_1,b_1)}  \bot \mathrm{Vect(a_0,b_0)}\\
       & a & \longmapsto & f(a)+(q(a)+q(f(a))) a_1+a_0\\
       & b & \longmapsto & f(b)+(q(b)+q(f(b))) a_1+ (1+B(f(a), f(b))) b_0
\end{array}$$
preserves the quadratic form. Consequently, the morphism:
$$T=[ V \xrightarrow{g_1} W \bot H_1 \bot H_0 \hookleftarrow W]$$
is a morphism of $\Tq$ such that $\epsilon(T)=f$.

Let $V_n$ be a non-degenerate quadratic space of dimension $2n$, $\{a_1,
b_1, \ldots, a_n, b_n \}$ be a symplectic basis of $V_n$ and $f_n: V_n
\rightarrow W$ be a linear map. By induction, there exists a  map :

$$ \begin{array}{cccc}
g_n: & V_n& \rightarrow & W \bot Y\\
       & a_1 & \longmapsto & f(a_1)+y_1\\
       & b_1 & \longmapsto & f(b_1)+z_1\\
       & \ldots& & \ldots\\
       & a_n & \longmapsto & f(a_n)+y_n\\
       & b_n & \longmapsto & f(b_n)+z_n\\
\end{array}$$
where $y_i$ and $z_i$, for all integers $i$ between $1$ and $n$, are
elements of $Y$, which preserves the quadratic form and such that: 
$$\epsilon([V_n \xrightarrow{g_n} W \bot Y \hookleftarrow W])=f_n.$$

Let $V_{n+1}$ be a non-degenerate quadratic space of dimension
$2(n+1)$,\\
$\{a_1,
b_1, \ldots, a_n, b_n, a_{n+1}, b_{n+1} \}$ a symplectic basis of
$V_{n+1}$ and $f_{n+1}: V_{n+1} \rightarrow W$ a linear map. To define the
map $g_{n+1}$, we will consider the restriction of $f_{n+1}$ over
$V_n$ and extend the map $g_n$ given by the inductive assumption. For
that, we need the following space: $E \simeq W \bot Y \bot H_0^{\bot n} \bot H_0^{\bot
  n} \bot H_1 \bot H_0$, for which we specify the notations for a basis:
$$ E \simeq  W \bot Y \bot ( \bot_{i=1}^n \mathrm{Vect}(a_0^i, b_0^i)) \bot ( \bot_{i=1}^n \mathrm{Vect}(A_0^i, B_0^i)) 
  \bot \mathrm{Vect}(A_1, B_1) \bot \mathrm{Vect}(C_0, D_0).$$
We verify that the following map: 
$$ \begin{array}{ccl}
 V& \xrightarrow{g_{n+1}} & W \bot Y \bot H_0^{\bot n} \bot H_0^{\bot n} \bot H_1 \bot H_0\\
        a_1 & \longmapsto & f(a_1)+y_1+a_0^1\\
        b_1 & \longmapsto & f(b_1)+z_1+A_0^1\\
        \ldots& & \ldots\\
        a_i & \longmapsto & f(a_i)+y_i+a_0^i\\
        b_i & \longmapsto & f(b_i)+z_i+A_0^i\\
        \ldots& & \ldots\\
        a_n & \longmapsto & f(a_n)+y_n+a_0^n\\
        b_n & \longmapsto & f(b_n)+z_n+A_0^n\\
        a_{n+1} & \longmapsto &
       f(a_{n+1})+(q(a_{n+1})+q(f(a_{n+1})))A_1+C_0+\sum_{i=1}^n
       B(f(a_i),f(a_{n+1})) b_0^i\\
            &          &+\sum_{i=1}^n
       B(f(b_i),f(a_{n+1})) B_0^i\\
        b_{n+1} & \longmapsto &
       f(b_{n+1})+(q(b_{n+1})+q(f(b_{n+1})))A_1+(1+B(f(a_{n+1}),
       f(b_{n+1})))D_0\\   
            &          &+\sum_{i=1}^n
       B(f(a_i), f(b_{n+1})) b_0^i+\sum_{i=1}^n
       B(f(b_i),f(b_{n+1})) B_0^i       
\end{array}$$
preserves the quadratic form and we have 
 $$\epsilon([V \xrightarrow{g_{n+1}} W \bot Y \bot H_0^{\bot n} \bot
 H_0^{\bot n} \bot H_1 \bot H_0\hookleftarrow W])=f$$
which completes the inductive step.
\end{proof}

\subsection{Proof of theorem \ref{3.1}}
Theorem \ref{3.1} is a consequence of a general result about functor
categories which we recall in Proposition \ref{A2} of the
appendix. But, as the functor $\epsilon$ is not essentially
surjective, we can not apply directly the proposition to the category $\F$. Consequently, we introduce a
category $\F'$ equivalent to $\F$.
\begin{defi}
The category $\E^{f-(even)}$ is the full subcategory of $\E^f$ having as
objects the $\FF$-vector spaces of even dimension.
\end{defi}
\begin{nota} \label{}
We denote by $\F'$ the category of functors from $\E^{f-(even)}$ to $\E$.
\end{nota}
We have the following result:
\begin{prop} \label{3.8}
The categories $\F$ and $\F'$ are equivalent.
\end{prop}
The proof of this proposition relies on the following standard lemma.
\begin{lm} \cite{KuhnII}
Let $n>0$ be a natural integer. Any idempotent linear map $e_{2n-1}: \FF^{2n}
\rightarrow \FF^{2n}$ of rank $2n-1$, verifies:
$$ P^{\F}_{\FF^{2n-1}} \simeq P^{\F}_{\FF^{2n}} \cdot e_{2n-1}$$
where $P^{\F}_{\FF^{n}}(-)=\FF[\mathrm{Hom}_{\E^f}(\FF^{n}, -)]$ is the
standard projective object of $\F$ given by the Yoneda lemma.
\end{lm}

\begin{proof}[Proof of proposition \ref{3.8}]
Let $V$ be an object of $\E^f$. 

If the dimension of $V$ is even, $V$ is an object of $\E^{f-(even)}$.

If the dimension of $V$ is odd $\mathrm{dim}(V)=2n-1$, by the previous
lemma, $P^{\F}_{\FF^{2n-1}}$ is a direct summand of
$P^{\F}_{\FF^{2n}}$.

As the category $\E^{f-(even)}$ is a full
subcategory of $\E^f$, we obtain by proposition \ref{A3} of the
appendix, the theorem.
\end{proof}

\begin{proof}[Proof of theorem \ref{3.1}]
The functor $\iota$ of theorem \ref{3.1} is, by definition, the
precomposition functor by the functor $\epsilon: \Eq \rightarrow
\E^f$. The two first points of the theorem are clear. As the objects of $\Eq$ are spaces of even dimension, the functor
$\epsilon:\Eq \rightarrow \E^{f}$ factorizes through the inclusion $
\E^{f-(even)} \hookrightarrow \E^{f}$. This induces a functor
$\epsilon': \Eq \rightarrow \E^{f-(even)}$ which is full and
essentially surjective. We deduce from proposition \ref{A2} that
the functor $- \circ
\epsilon':\F' \rightarrow \Fq$ is fully-faithful and, for a simple
object $S$ of $\F'$, $S \circ \epsilon'$ is simple in $\Fq$. The
theorem follows from proposition \ref{3.8}.
\end{proof}

\subsection{Duality}
In this section, we prove that the duality defined over $\Fq$ in
\ref{2.14} is an extension of the duality over $\F$, $D: \F^{op}
\rightarrow \F$, given by $DF=-^* \circ F \circ (-^*)^{op}$ for $F$ an object
of $\F$ and $-^*$ the duality functor from $\E^{op}$ to $\E$.
\begin{prop}
We have the following commutative diagram, up to natural isomorphism
$$\xymatrix{
   \F^{op} \ar[r]^{\iota^{op}} \ar[d]_{D} & \Fq^{op} \ar[d]^{D} \\
   \F \ar[r]_{\iota} & \Fq. }$$
\end{prop}
\begin{proof}
This relies on the following commutative diagram:
$$\xymatrix{
   \Tq^{op} \ar[r]^{\epsilon^{op}} \ar[d]_{tr} & (\E^f)^{op} \ar[d]^{-^*} \\
   \Tq \ar[r]_{\epsilon} & \E^f. }$$
The commutativity is a consequence of classical results of linear
algebra about the duality of vector spaces and the fact that a
non-singular bilinear form on a vector space $V$ determines a
privileged isomorphism between $V$ and $V^*$.
\end{proof}

\section{The category $\Gq$ }
In this section, we define a subcategory $\Gq$ of $\Fq$ which is, by
theorem \ref{4.7}, an abelian symmetric monoidal category with enough
projective objects. The
category $\Gq$ is related to $\Fq$ by the following theorem.

\begin{thm} \label{4.1}
There is a functor $$\kappa: \Gq \rightarrow \Fq$$ which satisfies the
following properties:
\begin{enumerate}
\item
$\kappa$ is exact;
\item
$\kappa$ preserves tensor products;
\item
$\kappa$ is fully-faithful;
\item
if $S$ is a simple object of $\Gq$, $\kappa(S)$ is a simple object of $\Fq$.
\end{enumerate} 
\end{thm}
We obtain a classification of the simple objects of $\Gq$ from the following theorem.
\begin{thm} \label{4.2}
There is a natural equivalence 
$$\Gq \simeq \prod_{V \in \mathcal{S}} \FF[O(V)]-mod$$
where $\mathcal{S}$ is a set of representatives of isometry classes of
objects of $\Eqd$.
\end{thm}
\subsection{Definition of the category $\Gq$}
\subsubsection{The category $\mathrm{Sp}(\Eqd)$}
Example \ref{2.3} implies that the category $\mathrm{Sp}(\Eqd)$
is defined. In this section, we give some properties of this category
which are similar to those given for the category $\Tq$.

\begin{defi}[Transposition] \label{4.3}
The transposition functor, $tr: \mathrm{Sp}(\Eqd)^{op}
\rightarrow \mathrm{Sp}(\Eqd)$, is defined on objects by $tr(V)=V$ and on
 morphisms by:
$$\mathrm{tr}(f)=\mathrm{tr} (\lbrack  V \leftarrow X \rightarrow W
\rbrack)= \lbrack  W \leftarrow X \rightarrow V \rbrack,$$
for $f$ an element of $\mathrm{Hom}_{\mathrm{Sp}(\Eqd)^{op}}(W,V)$.
\end{defi}

\begin{prop}[Orthogonal sum] \label{4.5}
There exists a bifunctor $\bot: \mathrm{Sp}(\Eqd) \times \mathrm{Sp}(\Eqd) \rightarrow \Sq$,
called the orthogonal sum, defined on objects by:
$$\bot (V,W)=V \bot W$$
and on morphisms by:
$$\bot (\lbrack V
\xleftarrow{f} X \xrightarrow{g} W \rbrack, \lbrack V'
\xleftarrow{f'} X' \xrightarrow{g'} W' \rbrack)=\lbrack V \bot V'
\xleftarrow{f \bot f'} X \bot X' \xrightarrow{g \bot g'} W \bot W'
\rbrack.$$
This bifunctor gives $\mathrm{Sp}(\Eqd)$ the structure of a symmetric
monoidal category, with unit $\{ 0 \}$.
\end{prop}

\subsubsection{The category $\Gq$}
As the category $\mathrm{Sp}(\Eqd)$ is essentially small, we can give
the following definition.
\begin{defi}
The category $\Gq$ is the category of functors from $\mathrm{Sp}(\Eqd)$ to $\E$.
\end{defi}
\begin{rem}
The category $\Gq$ is equivalent to the category
of Mackey functors from $\Eqd$ to $\E$, by the paper \cite{lindner}. 
\end{rem}
As for the category $\Fq$, we obtain the following
theorem.
\begin{thm} \label{4.7}
\begin{enumerate}
\item
The category $\Gq$ is abelian.
\item
The tensor product of vector spaces induces a structure of symmetric
monoidal category on $\Gq$.
\item
For an object $V$ of $\ \mathrm{Sp}(\Eqd)$, the functor $Q_V=\FF \lbrack
\mathrm{Hom}_{\mathrm{Sp}(\Eqd)}(V,-) \rbrack $ is a projective object
and there is a natural isomorphism:
$$\mathrm{Hom}_{\Gq}(Q_V, F) \simeq F(V)$$
for all objects $F$ of $\Gq$.

The set of functors $\{ Q_V | V \in \mathcal{S} \}$ is a set of
projective generators of $\Gq$, where $\mathcal{S}$ is a set of
representatives of isometry classes of degenerate quadratic spaces. In particular, the category $\Gq$ has enough projective objects.
\end{enumerate}
\end{thm}

\begin{defi}
The duality functor of $\Gq$ is the functor $D: \Gq^{op} \rightarrow
\Gq$ given by:
$$DF=-^* \circ F \circ tr^{op}$$
for $F$ an object of $\Gq$, $-^*$ the duality functor from $\E^{op}$ to
$\E$ and $tr$ is the transposition functor of $\mathrm{Sp}(\Eqd)$ defined in \ref{4.3}.
\end{defi}

The following proposition summarizes the basic properties of the duality
functor $D$.

\begin{prop} \label{4.9}
\begin{enumerate}
\item
The functor $D$ is exact.
\item
The functor $D$ is right adjoint to the functor $D^{op}$, i.e. we have
a natural isomorphism:
$$\mathrm{Hom}_{\Gq}(F, DG) \simeq \mathrm{Hom}_{\Gq^{op}}(D^{op} F, G)\simeq \mathrm{Hom}_{\Gq}(G, DF).$$
\item
For $F$ an object of $\Gq$ with values in finite dimensional vector
spaces, the unit of the adjunction between $\Gq$ and $\Gq^{op}$: $F
\rightarrow DD^{op}F$ is an isomorphism.

\end{enumerate}
\end{prop}
A straightforward consequence of the first point of the last
proposition and theorem \ref{4.7} is:
\begin{cor}
The category $\Gq$ has enough injective objects.
\end{cor}
\subsection{An equivalent definition of $\Gq$ }
In order to apply proposition \ref{A2} in the appendix to prove
 theorem \ref{4.1}, we will use the same strategy as for
theorem \ref{3.1}. In
other words, we will introduce a category equivalent to $\Gq$ such
that $\kappa$ will be the
precomposition functor by a full and essentially surjective functor.
First, we give the following definition.
\begin{defi}
The category $\Sq$ is the full subcategory of $\mathrm{Sp}(\Eqd)$
having as objects the non-degenerate quadratic spaces.
\end{defi}
\begin{rem}
A morphism of $\Sq$ is represented by a diagram $V \xleftarrow{f} D
\xrightarrow{g} W$ where $V$ and $W$ are non-degenerate quadratic
spaces and $D$ is a possibly degenerate quadratic space.
\end{rem}
The transposition functor and the orthogonal sum defined in the
previous section for the category $\mathrm{Sp}(\Eqd)$ induce, by
restriction, a transposition functor and an orthogonal sum for the category $\Sq$.

As the category $\Sq$ is, by definition, a full subcategory of the category
$\mathrm{Sp}(\Eqd)$, we have the existence of a functor $\lambda'$
from $\Gq$ to $\mathrm{Func}(\Sq, \E)$ induced by the inclusion $\Sq
\hookrightarrow \mathrm{Sp}(\Eqd)$. The aim of this section is to
show that $\lambda'$ is an isomorphism.

\begin{thm} \label{4.13}
There exists a natural isomorphism:
$$\Gq \simeq \mathrm{Func}(\Sq, \E)$$
where $\mathrm{Func}(\Sq, \E)$ is the category of functors from $\Sq$
to $\E$.
\end{thm}

To prove the theorem, we require some results about the idempotents
of the category $\mathrm{Sp}(\Eqd)$.

\subsubsection{The idempotents of $\mathrm{Sp}(\Eqd)$}
We begin this section with the following notation.
\begin{nota} 
Let $V$ be an object of $\Eqd$, $\alpha$ and $\beta$ be morphisms of
$\mathrm{Hom}_{\Eqd}(D,V)$. We denote by $f_{\alpha, \beta}$ the
following morphism of $\mathrm{Hom}_{\mathrm{Sp}(\Eqd)}(V,V)$:
$$[V \xleftarrow{\alpha} D \xrightarrow{\beta} V].$$
\end{nota} \label{lm}
\begin{prop} \label{4.15}
\begin{enumerate}
\item
An idempotent of $\mathrm{Hom}_{\mathrm{Sp}(\Eqd)}(V,V)$ is of the form 
$$e_\alpha=[V \xleftarrow{\alpha} D \xrightarrow{\alpha} V]$$
where $\alpha$ is an element of
$\mathrm{Hom}_{\Eqd}(D,V)$, for some $D$.
\item
For $\alpha$ and $\beta$ two morphisms of $\Eqd$ with range $V$,
the idempotents $e_{\alpha}$ and $e_{\beta}$ commute.
\item
For $\alpha$ and $\beta$ two morphisms of $\Eqd$ with range $V$,
the elements $1+e_{\alpha}$ and $1+e_{\beta}$ are idempotents of
$\FF[\mathrm{Hom}_{\mathrm{Sp}(\Eqd)}(V,V)]$ which commute. 
\end{enumerate}
\end{prop}
\begin{proof}
By definition of $\mathrm{Sp}(\Eqd)$, $f_{\alpha, \beta} \circ f_{\alpha, \beta}=[V
\leftarrow D' \rightarrow V]$ where
$D'$ is the pullback in $\Eqd$ of the diagram $D \xrightarrow{\beta} V
\xleftarrow{\alpha} D$ and $f_{\alpha, \beta} \circ f_{\alpha,
  \beta}=f_{\alpha, \beta}$ if and only if there exists an isomorphism
$g: D \rightarrow D'$ such that the following diagram is commutative:
$$\xymatrix{
D \ar[rrrd]^{\beta} \ar[rd]^{\simeq}_{g} \ar[rddd]_{\alpha}\\
    &  D' \ar[r]^{g^{-1}} \ar[d]^{g^{-1}}& D \ar[r]_{\beta}
    \ar[d]^{\alpha}& V\\
    & D \ar[r]_{\beta} \ar[d]^{\alpha} & V\\
    &V
}$$
This implies that $\alpha =\beta$. Consequently
$f_{\alpha, \beta} \circ f_{\alpha, \beta}=f_{\alpha, \beta}$ if and
only if $f_{\alpha, \beta}= e_{\alpha}$.

The second point is straightforward and the last point is a direct
consequence of the second one by a standard result about
idempotents (noting that $-1=1$ in $\FF$).
\end{proof}

\subsubsection{Proof of theorem \ref{4.13}}

The proof of theorem \ref{4.13} relies on the following crucial lemma.

\begin{lm} \label{4.16}
For $V$ an object of $\mathrm{Span}(\Eqd)$ and $\alpha: A
\hookrightarrow V$ a subobject of $V$ in $\Eqd$, 
$$Q_V \cdot e_{\alpha} \simeq Q_A.$$
\end{lm}
\begin{proof}
Let $\alpha_*: Q_A \rightarrow Q_V$ (respectively $\alpha^*: Q_V
\rightarrow Q_A$) be the morphism of $\Gq$ which corresponds
by the Yoneda lemma to the element $[V \xleftarrow{\alpha} A
\xrightarrow{\mathrm{Id}} A]$ of $Q_V(A)$ (respectively $[A \xleftarrow{\mathrm{Id}} A
\xrightarrow{\alpha} V]$ of $Q_A(V)$). 

As 
$[V \xleftarrow{\alpha} A
\xrightarrow{\mathrm{Id}} A] \circ [A \xleftarrow{\mathrm{Id}} A
\xrightarrow{\alpha} V]=\mathrm{Id}$
we have 
$\alpha^* \circ \alpha_* =\mathrm{Id}$ and as $[A \xleftarrow{\mathrm{Id}} A
\xrightarrow{\alpha} V] \circ [V \xleftarrow{\alpha} A
\xrightarrow{\mathrm{Id}} A] =e_{\alpha}$ we have $\alpha_* \circ
\alpha^*=\cdot e_{\alpha}$.

\end{proof}

\begin{proof}[Proof of theorem \ref{4.13}]
Let $A$ be an object of $\mathrm{Sp}(\Eqd)$, there exists an object
$V$ of $\Sq$ such that $A$ is a subobject of $V$. By the
previous lemma, we deduce that $Q_A$ is a direct summand of $Q_V$. As 
the category $\Sq$ is a full subcategory of $\mathrm{Sp}(\Eqd)$ we
obtain, by proposition \ref{A3} of the appendix, the theorem.\end{proof}

\subsection{Relation between $\Gq$ and $\Fq$}
The main result of this section is theorem \ref{4.1}, which gives the
existence of an exact, fully-faithful functor $\kappa: \Gq \rightarrow
\Fq$ which preserves the simple objects. To define
the functor $\kappa$ of this theorem we need to define and study the
functor $\sigma: \Tq \rightarrow \Sq$.
\subsubsection{Definition of  $\sigma: \Tq \rightarrow \Sq$}
\begin{prop} \label{4.17}
There exists a monoidal functor $\sigma: \Tq \rightarrow \Sq$ defined by $\sigma(V)=V$
 and
$$\sigma (\lbrack V \rightarrow X \leftarrow W \rbrack)=\lbrack V
\leftarrow V \underset{X}{\times} W \rightarrow W\rbrack$$
where $V \underset{X}{\times} W$ is the pullback in $\Eqd$.
\end{prop}
The proof of this proposition relies on the following important lemma.
\begin{lm} \label{4.18}
For $V$, $W$ and $X$ objects of $\Eq$ and $V \xrightarrow{f} W$, $V
\xrightarrow{g} X$ morphisms of $\Eq$, we have:
$$X  \underset{X \underset{V}{\bot} W}{\times} W \simeq V$$
where $-  \underset{A}{\times} -$ is the pullback over $A$ in
$\mathcal{E}_{q}^{\mathrm{deg}}$ and $-  \underset{B}{\bot} -$ is the
pseudo push-out over $B$ in $\Eq$.
\end{lm}
\begin{proof}
It is a straightforward consequence of the definitions of the pseudo
push-out and the pullback.
\end{proof}
\begin{rem}
The result of the previous lemma explains why we have imposed that the
morphisms of the category $\Eqd$ are monomorphisms.
\end{rem}
\begin{proof}[Proof of proposition \ref{4.17}]
The functor $\sigma$ is well defined on
the classes of morphisms of $\Tq$. To show that $\sigma$ is a
functor, we verify that 
$$\sigma([V \xrightarrow{\mathrm{Id}} V \xleftarrow{\mathrm{Id}}
V])=\mathrm{Id}$$
and, that $\sigma$ respects composition, which is a direct consequence of the
above lemma. Finally, the following consequence of the
definition of the pullback and of the orthogonal sum:\\
for $T_1 =\lbrack V_1 \rightarrow X_1 \leftarrow W_1
   \rbrack$ and  $T_2 =\lbrack V_2 \rightarrow X_2 \leftarrow W_2
   \rbrack$, we have:

$$(V_1 \underset{X_1}{\times} W_1) \bot (V_2 \underset{X_2}{\times}
W_2) \simeq (V_1 \bot V_2) \underset{X_1 \bot X_2}{\times} (W_1 \bot
W_2)$$
implies that $\sigma$ preserves the
monoidal structures.
\end{proof}
\subsubsection{The fullness of $\sigma$}
The aim of this section is to prove the following proposition.
\begin{prop} \label{4.20}
The functor $\sigma$ is full.
\end{prop}
The proof of this proposition relies on the following technical lemma.
\begin{lm} \label{4.21}
Let $D$ be a degenerate quadratic space of dimension $r$ which has the
following decomposition
$$D=(x_1, \epsilon_1) \bot \ldots \bot (x_r,\epsilon_r)$$
where $\epsilon_i \in \{0,1 \} $, $H$ a non-degenerate quadratic space
and $f$ an element of \\
$\mathrm{Hom}_{\Eqd}(D,H)$. Then there exists
elements $k_1, \ldots ,k_r$ in $H$ and a non-degenerate quadratic
space $H'$ such that  
$$H=\mathrm{Vect}(f(x_1),k_1)\bot \ldots \bot
\mathrm{Vect}(f(x_r),k_r)\bot H'$$ and 
$B(f(x_i),k_i)=1$ where  $B$ is the underlying bilinear form.
\end{lm}
\begin{proof}
We prove this lemma by induction on the dimension $r$ of the space
$D$.

For $r=1$, we have $f: (x_1,\epsilon_1) \rightarrow H$ and
$f(x_1)=h_1$. As $H$ is, by hypothesis, non-degenerate there exists an
element $k_1$ in $H$ such that $B(h_1,
k_1)=1$. Then, the space  $K=\mathrm{Vect}(h_1,k_1)$ is a
non-degenerate subspace of $H$, so we have $H=K \bot K^{\bot}$, with
$K^{\bot}$ non-degenerate.

Suppose that the result is true for $r=n$. Let $(x_1, \ldots, x_n,
x_{n+1})$ be linearly independent vectors and $f: (x_1,
\epsilon_1) \bot \ldots \bot (x_n,\epsilon_n) \bot (x_{n+1},
\epsilon_{n+1}) \rightarrow H$. By restriction, we have 
$$ f \circ i:(x_1,
\epsilon_1) \bot \ldots \bot (x_n,\epsilon_n) \rightarrow H$$
 and by the inductive assumption, we have the existence of $k_1,
 \ldots, k_n$ in $H$ such that:
$$H=\mathrm{Vect}(f(x_1),k_1) \bot \ldots \mathrm{Vect} (f(x_n),k_n) \bot H'.$$

We decompose $f(x_{n+1})$ over this basis to obtain the following decomposition
$$f(x_{n+1})=\Sigma_{i=1}^n (\alpha_i f(x_i)+\beta_i k_i)+h'$$
where $\alpha_i$ and $\beta_i$ are elements of $\FF$ and $h'$ is
an element of $H'$.
As $f$ preserves the quadratic form and, consequently, the underlying
bilinear form, we have $\beta_i=B(f(x_{n+1}),
f(x_i))=B(x_{n+1},x_i)=0$ for all $i$. Hence
\begin{equation} \label{5}
f(x_{n+1})=\Sigma_{i=1}^n \alpha_i f(x_i)+h'.
\end{equation}

As the vectors $(f(x_1), \ldots, f(x_n),
f(x_{n+1}))$ are linearly independent, by the injectivity of $f$, we have $h' \ne 0$
and, as $H'$ is non-degenerate, there exists an element $k'$ in $H'$
such that $B(h',k')=1$. We deduce the following decomposition $H'=\mathrm{Vect}(h',k')
\bot H''$.

In the equality (\ref{5}), after reordering, we can suppose that 
$$ \alpha_1=\ldots =\alpha_p=1$$
and
$$ \alpha_{p+1}=\ldots =\alpha_n=0.$$

Then we have, 
$$ B(f(x_{n+1}),k_i)=1 \ \mathrm{pour\ } i=1, \ldots, p$$
and
$$ B(f(x_{n+1}),k_i)=0 \ \mathrm{pour\ } i=p+1, \ldots, n.$$

Consequently, by the following decomposition of $H$
$$H=\bot_{i=1}^{p} \mathrm{Vect}(f(x_i),k_i +k') \bot_{j=p+1}^{n}
\mathrm{Vect}(f(x_{j}),k_{j}) \bot \mathrm{Vect}(f(x_{n+1}),k') \bot H''$$
we obtain the result.
\end{proof}

\begin{proof}[Proof of proposition \ref{4.20}]
We prove that, for a morphism  $S =\lbrack V
   \xleftarrow{f} D  \xrightarrow{g} W  \rbrack$ of  $\Sq$, there
   exists a morphism $T$ in $\mathrm{Hom}_{\Tq}(V,W)$ such that 
$\sigma (T)=S$.

First we decompose the morphism $S$ as an orthogonal sum of more
simple morphisms. By theorem \ref{1.9} we have $D \simeq H \bot \mathrm{Rad}(V)$
such that $H$ is non-degenerate and  $\mathrm{Rad}(V)$ is isometric to
either $(x,0)^{\bot r}$ or $(x,1)^{\bot r}$, where $r$ is the
dimension of $\mathrm{Rad}(V)$. By the previous lemma, we can
decompose the morphisms $f$ and $g$ in the following form:
$$ f:\ H \bot  \mathrm{Rad}(D) \rightarrow H \bot D' \bot V' \simeq
V$$
and
$$ g:\ H  \bot \mathrm{Rad}(D) \rightarrow H \bot D'' \bot W' \simeq
W$$
where $D'$ (respectively $D''$) is one of the non-degenerate spaces
constructed in lemma \ref{4.21} and $V'$ (respectively $W'$) is the
orthogonal of $H \bot D'$ (respectively $H \bot D''$). We deduce that
:
$$S= \lbrack H \leftarrow H  \rightarrow H  \rbrack \bot \lbrack V'
\leftarrow 0  \rightarrow W' \rbrack \bot \lbrack D'
\leftarrow \mathrm{Rad}(D) \rightarrow D'' \rbrack. $$

Furthermore, again by theorem \ref{1.9}
$$\lbrack D'
\leftarrow \mathrm{Rad}(D) \rightarrow D'' \rbrack =\bot_{i}\lbrack
V_i \leftarrow (x_i, \epsilon_i) \rightarrow W_i \rbrack   $$
where, by lemma \ref{4.21}, $V_i$ and $W_i$ are spaces of dimension
two, therefore:
$$S= \lbrack H \leftarrow H  \rightarrow H  \rbrack \bot \lbrack V'
\leftarrow 0  \rightarrow W' \rbrack \bot \lbrack V_1
\leftarrow (x_1,\epsilon_1) \rightarrow W_1 \rbrack \bot \ldots \bot \lbrack V_r
\leftarrow (x_r,\epsilon_r) \rightarrow W_r\rbrack. $$

According to proposition \ref{4.17}, it is enough to prove that for each 
morphism $S_{\alpha}$, which appears as a factor in the previous decomposition
of $S$, there exists a morphism of $\Tq$, $T_{\alpha}$,
such that $\sigma(T_{\alpha})=S_{\alpha}$.

Obviously, we have
$$\sigma (\lbrack H \rightarrow H  \leftarrow H  \rbrack) = \lbrack H
\leftarrow H  \rightarrow H  \rbrack$$
and 
$$\sigma (\lbrack V \rightarrow V \bot W  \leftarrow W  \rbrack) = \lbrack V
\leftarrow 0  \rightarrow W  \rbrack.$$

For the morphisms $\lbrack V
\leftarrow (x, \epsilon) \rightarrow W \rbrack$, we have to consider
several cases.

In the case $\epsilon =0$, as all the non-zero element $x$ of $H_1$
verify $q(x)=1$, we have
$\mathrm{Hom}_{\Eq^{\mathrm{deg}}}((x,0),H_1)=
\emptyset$. Consequently, we have to consider only the following morphism:
$$S_1=\lbrack H_0
\xleftarrow{f} (x,0)  \xrightarrow{g} H_0  \rbrack.$$

After composing by an element of $O^+_2=O(H_0)$, we can suppose that
$f(x)=g(x)=a_0$. The morphism 
$$T_1=\lbrack H_0 \xrightarrow{f'} H_0 \bot H_0
\xleftarrow{g'} H_0  \rbrack$$
where, $f'$ is defined by:
$$f'(a_0)=a_0\ \mathrm{and}\ f'(b_0)=b_0+a'_0$$
and $g'$ is defined by:
$$g'(a_0)=a_0\ \mathrm{and}\ g'(b_0)=b_0+b'_0$$
satisfies $\sigma(T_1)=S_1$.\\
$\ $

In the case $\epsilon =1$, we have to consider three kinds of morphisms.

For $$S_2=\lbrack H_0
\xleftarrow{f} (x,1)  \xrightarrow{g} H_0  \rbrack,$$
as only the element $a_0+b_0$ of $H_0$ verify $q(a_0+b_0)=1$, we have
$f(x)=g(x)=a_0+b_0$. The morphism 
$$T_2=\lbrack H_0 \xrightarrow{f'} H_0 \bot H_0
\xleftarrow{g'} H_0  \rbrack$$
where, $f'$ is defined by:
$$f'(a_0)=a_0+a'_0\ \mathrm{and}\ f'(b_0)=b_0+a'_0$$
and $g'$ is defined by:
$$g'(a_0)=a_0+b'_0\ \mathrm{and}\ g'(b_0)=b_0+b'_0$$
satisfies $\sigma(T_2)=S_2$.\\
$\ $

For $$S_3=\lbrack H_1
\xleftarrow{f} (x,1)  \xrightarrow{g} H_1  \rbrack,$$
after composing by an element of $O^-_2=O(H_1)$, we can suppose that
$f(x)=g(x)=a_1$. The morphism  
$$T_3=\lbrack H_1 \xrightarrow{f'} H_1 \bot H_0
\xleftarrow{g'} H_1  \rbrack$$
where, $f'$ is defined by:
$$f'(a_1)=a_1\ \mathrm{and}\ f'(b_1)=b_1+b_0$$
and $g'$ is defined by:
$$g'(a_1)=a_1\ \mathrm{and}\ g'(b_1)=b_1+a_0$$
satisfies $\sigma(T_3)=C_3$.\\
$\ $

For $$S_4=\lbrack H_0
\xleftarrow{f} (x,1)  \xrightarrow{g} H_1  \rbrack,$$
we have $f(x)=a_0+b_0$ and, after composing by an element of $O^-_2=O(H_1)$, we can suppose that $g(x)=a_1$.  The morphism   
$$T_4=\lbrack H_0 \xrightarrow{f'} H_1 \bot H_0
\xleftarrow{g'} H_1  \rbrack$$
where, $f'$ is defined by:
$$f'(a_0)=a_1+b_1+a_0+b_0\ \mathrm{et}\ f'(b_0)=b_1+a_0+b_0$$
and $g'$ is defined by:
$$g'(a_1)=a_1\ \mathrm{et}\ g'(b_1)=b_1$$
satisfies $\sigma(T_4)=C_4$.

The final possibility results from the previous one by transposition. 
\end{proof}

\subsubsection{Proof of theorem \ref{4.1}}

The functor $\kappa$ of the theorem \ref{4.1} is, by definition, the
precomposition functor by $\sigma: \Tq \rightarrow \Sq$. The two first
point of the theorem are clear. By the
proposition \ref{4.17} the functor $\sigma$ is full and, by definition,
it is essentially surjective. So, the two last points of the theorem
\ref{4.1} are direct consequences of the proposition \ref{A2} given in
the appendix.

\subsubsection{Duality}
\begin{prop}
We have the following commutative diagram, up to natural isomorphism
$$\xymatrix{
   \Gq^{op} \ar[r]^{\kappa^{op}} \ar[d]_{D} & \Fq^{op} \ar[d]^{D} \\
   \Gq \ar[r]_{\kappa} & \Fq. }$$
\end{prop}

\begin{proof}
This relies on the following commutative diagram:
$$\xymatrix{
   \Tq^{op} \ar[r]^-{\sigma^{op}} \ar[d]_{tr} & (\mathrm{Sp}(\Eqd))^{op} \ar[d]^{tr} \\
   \Tq \ar[r]_-{\sigma} & \mathrm{Sp}(\Eqd) }$$
which is a direct consequence of the definitions.
\end{proof}

\subsection{The isotropic functors}
In this section we are interested in an important family of functors of
$\Gq$ named the isotropic functors; the choice of terminology will
be explained below. After giving the definition of these functors we prove that they
are self-dual. We will show in the following sections that these
functors give rise to a family of projective generators of the
category $\Gq$.
\subsubsection{Definition}
Let $(\mathrm{Id}_V)^*$ be the element of 
$$DQ_V(V)=Q_V(V)^*=\mathrm{Hom}(\FF[\mathrm{End}_{\mathrm{Sp}(\Eqd)}(V)],\FF)$$
defined by:
$$(\mathrm{Id}_V)^*([\mathrm{Id}_V])=1 \quad \mathrm{and} \quad
(\mathrm{Id}_V)^*([f])=0 \mathrm{\ for\ all\ } f \ne \mathrm{Id}_V.$$
We denote by $a_V: Q_V \rightarrow DQ_V$ the morphism of $\Gq$ which
corresponds by the Yoneda lemma to the element $(\mathrm{Id}_V)^*$ of
$DQ_V(V)$.
\begin{defi}
The isotropic functor $iso_V: \mathrm{Sp}(\Eqd)\rightarrow \E$ of
$\Gq$ is the image of $Q_V$ by the morphism $a_V$.
\end{defi} 
\begin{nota} \label{4.24}
We denote by $K_V$ the kernel of $a_V$. We have the following short
exact sequence:
$$0 \rightarrow K_V \rightarrow Q_V \xrightarrow{a_V} iso_V \rightarrow 0.$$
\end{nota}
In the following lemma, we give, an explicit description of the vector
spaces $K_V(W)$ and $iso_V(W)$, which are elementary consequences of
the definition.
\begin{lm} \label{4.25}
For an object $W$ of $\mathrm{Sp}(\Eqd)$, we have that:
\begin{itemize}
\item
$K_V(W)$ is the subvector space of $Q_V(W)$ generated by the elements
$[ V \leftarrow H \rightarrow W]$ where $H \not \simeq V$;
\item
as a vector space, $iso_V(W)$ is isomorphic to the subspace of $Q_V(W)$ generated
by the elements $[ V \xleftarrow{\mathrm{Id}} V \rightarrow W]$. Consequently
$iso_V(W)$ has basis the set $\mathrm{Hom}_{\Eqd}(V,W)$. 
\end{itemize}
\end{lm}
\begin{rem}
Observe that the isomorphism given in the second point is not natural.
\end{rem}

\begin{rem}
The terminology ``isotropic functor'' was motivated by the first case
considered by the author. For the quadratic space $(x,0)$, we have:
$$iso_{(x,0)}(V) \simeq \FF[I_V]$$
where  $I_V=\{ v\in V \setminus \{0 \}|\ q(v)=0 \}$ is the isotropic
cone of the quadratic space $V$.
\end{rem}
\subsubsection{Self-duality}
\begin{rem}
For simplicity, we will denote in this section, $D^{op}$ by $D$.
\end{rem}
To begin we recall several definitions.
\begin{defi}
\begin{enumerate}
\item
A morphism $b: F \rightarrow DF$ is self-adjoint if $b=Db \circ
\eta_F$, where $\eta_F: F \rightarrow
D^2 F$ is the unit of the adjunction between $D$ and $D^{op}$.
\item
A functor $F$ is self-dual if there exists an isomorphism
$$\gamma: F \xrightarrow{\simeq} DF$$
which is self-adjoint.
\end{enumerate}
\end{defi}
The main result of this section is the following proposition.
\begin{prop} \label{4.28}
The isotropic functors of $\Fq$ are self-dual.
\end{prop}
The proof of this proposition relies on the following lemma.
\begin{lm}
Let $F$ be a functor which takes finite dimensional values and $a: F
\rightarrow DF$ a self-adjoint morphism, then
$im(a)$ is self-dual.
\end{lm}
\begin{rem}
This lemma and its proof are direct consequences of lemma $1.2.4$
in \cite{piriou}.
\end{rem}
\begin{proof}
The morphism $a$ admits the following factorisation by $im(a)$
$$\xymatrix{
a:\ F \ar@{->>}[r]^-{p}&  im(a) \ar@{^{(}->}[r]^{j} & DF.
}$$
We deduce the existence of $\tilde{a}$ which makes  the following
diagram commutative
$$\xymatrix{
0 \ar[r] & \mathrm{Ker}(a) \ar[r] \ar[d] & F \ar[r]^-{p}
\ar[d]^{\eta_F}&  im(a) \ar[r] \ar[d]^{\tilde{a}} & 0\\
0 \ar[r] & \mathrm{Ker}(Dj) \ar[r] & D^2F \ar[r]^-{Dj}& D (im(a)) \ar[r] & 0.
}$$
As $F$ has finite dimensional values, the unit of the adjunction $\eta_F: F
\rightarrow D^2F$ is an isomorphism. Consequently, by the following
commutative diagram 

$$\xymatrix{
F \ar@{->>}[r]^p \ar[d]^{\eta_F}_{\simeq} & im(a) \ar@{^{(}->}[r]^-{j}
\ar[d]^{\tilde{a}}&  DF  \ar@{=}[d]\\
D^2F \ar@{->>}[r]_-{Dj} & D(im(a)) \ar@{^{(}->}[r]_-{Dp}& DF }$$
we obtain that $\tilde{a}$ is an isomorphism.

 If we dualize the first commutative diagram, we obtain:
$$Dp \circ D\tilde{a} = D \eta_F \circ D^2 j.$$
So
$$Dp \circ D\tilde{a} \circ \eta_{im(a)} \circ p = D \eta_F \circ D^2 j \circ
\eta_{im(a)} \circ p= j \circ p= a.$$
We have also
$$Dp \circ \tilde{a} \circ p= Dp \circ Dj \circ \eta_F=Da \circ
\eta_F=a.$$
We deduce that $D\tilde{a} \circ \eta_{im(a)}=\tilde{a}$.
\end{proof}
\begin{proof}[Proof of proposition \ref{4.28}]
By the previous lemma, it is enough to prove that the morphism
$a_V: Q_V \rightarrow DQ_V$ is self-dual. By the Yoneda lemma we have
the following commutative diagram:
$$\xymatrix{
\mathrm{Hom}_{\Gq}(Q_V, DQ_W) \ar[r]^{f} \ar[d]_{\simeq} &
\mathrm{Hom}_{\Gq}(Q_W, DQ_V) \ar[d]^{\simeq}\\
\FF[\mathrm{Hom}_{\mathrm{Sp}(\Eqd)}(W,V)]^* \ar[r]_{\FF[tr]^*} &
\FF[\mathrm{Hom}_{\mathrm{Sp}(\Eqd)}(V,W)]^*}$$
where $tr: \mathrm{Sp}(\Eqd)^{op}
\rightarrow \mathrm{Sp}(\Eqd)$ is the transposition functor defined in
\ref{4.3} and $f$ is the natural isomorphism given in proposition
\ref{4.9}. By definition, $a_V: Q_V \rightarrow DQ_V$ corresponds, by
the Yoneda lemma, to the element $(\mathrm{Id}_V)^*$ of
$DQ_V(V)$. Since $tr(\mathrm{Id}_V)=\mathrm{Id}_V$, we deduce that
$a_V=Da_V \circ \eta_{Q_V}$.
\end{proof}

\subsection{Decomposition of the projective objects $Q_V$ of $\Gq$}
The aim of this section is to prove the following theorem.
\begin{thm} \label{4.31}
For $V$ an object of $\Eqd$, we have:
$$Q_V \simeq \bigoplus_{
\begin{array}{c}
A \in \mathcal{S}_V
\end{array}}
 iso_A$$
where $\mathcal{S}_V$ is the set of subobjects of $V$ in $\Eqd$,
represented by a morphism $\alpha: A \hookrightarrow V$.
\end{thm}
The proof of this theorem relies on the following proposition.
\begin{prop} \label{4.32}
For an object $V$ of $\ \Eqd$, the element $E_V$ of
$\FF[\mathrm{Hom}_{\mathrm{Sp}(\Eqd)}(V,V)]$ 
$$E_V = \prod_{
\begin{array}{c}
\alpha: A \hookrightarrow V\\
A \in \mathcal{S}_V \setminus V
\end{array}}
(1+e_{\alpha})$$
verifies:
\begin{enumerate}
\item
$E_V \cdot E_V=E_V$ 
\item
$Q_V \cdot E_V \simeq iso_V$. In particular $\mathrm{Hom}(iso_V,F)
\simeq F(E_V) \cdot F(V).$
\end{enumerate}
\end{prop}
\begin{proof}
The first point is a direct consequence of proposition \ref{4.15}
$(3)$ and $(1)$.

For the second point, we consider the following split short exact sequence:
$$0 \rightarrow Q_V \cdot (1+E_V) \rightarrow Q_V \rightarrow Q_V
\cdot E_V \rightarrow 0$$
and we recall that we have by notation \ref{4.24} the following short
exact sequence:
$$0 \rightarrow K_V \rightarrow Q_V \xrightarrow{a_V} iso_V
\rightarrow 0.$$
By the $5$-lemma, to obtain the result, it is sufficient, to prove
that:
$$Q_V \cdot (1+E_V) \simeq K_V.$$
By expanding $1+E_V$, we obtain that 
$$1+E_V=\prod_{
\begin{array}{c}
\gamma: A \hookrightarrow V\\
A \in \mathcal{R}_V
\end{array}}
e_{\gamma}$$
where $\mathcal{R}_V$ is a subset of $\mathcal{S}_V \setminus
V$. Consequently 
$$Q_V \cdot (1+E_V) \subset \sum_{
\begin{array}{c}
\gamma: A \hookrightarrow V\\
A \in \mathcal{S}_V \setminus V
\end{array}}
Q_V \cdot e_{\gamma}=K_V$$
where the last equality is a direct consequence of lemma \ref{4.25}.

On the other side, for $\gamma: A \hookrightarrow V$ with $A \in
\mathcal{S}_V \setminus V$, we have 
$$(Q_V \cdot e_{\gamma}) \cdot E_V=Q_V \cdot e_{\gamma} \cdot
(1+e_{\gamma}) \prod_{
\begin{array}{c}
\alpha: A \hookrightarrow V\\
\alpha \neq \gamma\\
A \in \mathcal{S}_V \setminus V
\end{array}}
(1+e_{\alpha}) =0$$
where the first equality follows from proposition \ref{4.15}
$(3)$. 

Consequently $Q_V \cdot e_{\gamma} \subset Q_V \cdot (1+E_V)$
and 
$$K_V = \sum_{
\begin{array}{c}
\gamma: A \hookrightarrow V\\
A \in \mathcal{S}_V \setminus V
\end{array}}
Q_V \cdot e_{\gamma} \subset Q_V \cdot (1+E_V).$$

We deduce that 
$$Q_V \cdot (1+E_V) \simeq K_V.$$

\end{proof}
An important consequence of the previous proposition is given in the
following corollary.
\begin{cor} \label{4.33}
For $V$ and $W$ objects of $\Eqd$, we have:
$$\mathrm{Hom}_{\Gq}(iso_V, iso_W)
\simeq \left\lbrace
\begin{array}{cc}
\FF[O(V)] &\mathrm{if\ } W \simeq V\\
0 & \mathrm{otherwise.}
\end{array}
\right.
$$
\end{cor}
\begin{proof}
By lemma \ref{4.25} an element of $iso_W(V)$ is represented by a
linear sum of 
$$[W \leftarrow W \xrightarrow{f} V].$$
\begin{itemize}
\item{If $W \not \simeq V$.}
We have
$$(1+e_f)[W \leftarrow W \xrightarrow{f} V]=[W \leftarrow W
\xrightarrow{f} V]+[W \leftarrow W \xrightarrow{f} V]=0$$
as the idempotents $(1+e_f)$ commute by proposition \ref{4.15}, we deduce that:
$$\mathrm{Hom}(iso_V, iso_W)=iso_W(E_V)\  iso_W(V)=0$$
where the first equality is given by proposition \ref{4.32} $(2)$.
\item{If $W \simeq V$.}
For $\alpha: A \hookrightarrow V$ with $A \in \mathcal{S}_V \setminus
V$ we have:
$$e_{\alpha} [W \leftarrow W \rightarrow V]=[W \leftarrow A \rightarrow
V] \in K_W(V).$$ 
So $iso_W(e_{\alpha})[W \leftarrow W \rightarrow V ]=0$ and
$$\mathrm{Hom}(iso_V, iso_W)=iso_W(E_V)\ iso_W(V)=iso_W(V) \simeq \mathrm{Hom}_{\Eqd}(W,V) \simeq
\FF[O(V)]$$
by lemma \ref{4.25}.
\end{itemize}
\end{proof}

\begin{prop} \label{4.34}
For $\alpha: A \hookrightarrow V$ a subobject of $V$, the idempotent
$E_{\alpha}$ defined by:
$$E_{\alpha}=e_{\alpha} \sum_{
\begin{array}{c}
\beta: B \hookrightarrow A\\
B \in \mathcal{S}_A \setminus A
\end{array}}
(1+e_{\beta})$$
verifies 
$$Q_V \cdot E_{\alpha} \simeq iso_A.$$
\end{prop}
\begin{proof}
By the proposition \ref{4.15}, $E_{\alpha}$ is clearly an
idempotent. The result follows from proposition
\ref{4.32} and lemma \ref{4.16}.
\end{proof}

\begin{proof}[Proof of theorem \ref{4.31}]
By the proof of \ref{4.32}, there exists an exact sequence
$$\bigoplus_{
\begin{array}{c}
\gamma: A \hookrightarrow V\\
A \in \mathcal{S}_V
\end{array}}
Q_V  \cdot e_{\gamma} \rightarrow Q_V \rightarrow iso_V \rightarrow 0$$
hence a complex 
$$\bigoplus_{
\begin{array}{c}
\gamma: A \hookrightarrow V\\
A \in \mathcal{S}_V
\end{array}}
iso_A  \rightarrow Q_V \rightarrow iso_V \rightarrow 0$$
by proposition \ref{4.34}. We deduce from proposition \ref{4.15} that
the idempotents given in the
proposition \ref{4.34} are orthogonal. Consequently the map 
$$\bigoplus_{
\begin{array}{c}
\gamma: A \hookrightarrow V\\
A \in \mathcal{S}_V
\end{array}}
iso_A  \rightarrow Q_V $$ is injective. Furthermore, by lemma
\ref{4.25}, for an object $X$ of $\mathrm{Sp}(\Eqd)$
$$Q_V(X) \simeq iso_V(X) \oplus \Big( \bigoplus_{
\begin{array}{c}
\alpha: W \hookrightarrow V\\
W \in \mathcal{S}_V \setminus V
\end{array}} iso_W(X) \Big).$$
We deduce that the previous complex is a short exact
sequence, which is
split by proposition \ref{4.32}.
\end{proof}
By theorem \ref{4.7} and the self-duality of the isotropic
functors given in proposition \ref{4.28}, we obtain the following corollary of
 theorem \ref{4.31}.
\begin{cor} \label{4.36}
The set of functors $\{ iso_V | V \in \mathcal{S} \}$ is a set of
projective generators (resp. injective cogenerators) of $\Gq$, where $\mathcal{S}$ is a set of
representatives of isometry classes of possibly degenerate quadratic spaces.
\end{cor}
\subsection{Proof of theorem \ref{4.2}}
In this section we prove the equivalence between the category $\Gq$
and the product of categories of modules over the orthogonal
groups.
\begin{proof}[Proof of theorem \ref{4.2}]
In \cite{Popescu} we have the following result (Corollary $6.4$, p. 103).

For any abelian category $\C$ the following assertions are equivalent:
\begin{enumerate}
\item
The category $\C$ has arbitrary direct sums and $\{ P_i \}_{i \in I}$
is a set of projective generators of finite type of $\C$.
\item
The category $\C$ is equivalent to the subcategory $\mathrm{Func}^{add}(\mathcal{P}^{op}, \mathcal{A}b)$  of\\
$\mathrm{Func}(\mathcal{P}^{op}, \mathcal{A}b)$ having as objects the functors satisfying
$F(f+g)=F(f)+F(g)$ where $f$ and $g$ are morphisms of
$\mathrm{Hom}_{\mathcal{P}^{op}}(V,W)$ and $\mathcal{P}$ is
the full subcategory of $\C$ having as objects $\{ P_i \ |\ i \in I \}$.
\end{enumerate}

Let $\C$ be the category $\Gq$. By corollary \ref{4.36}, the set of functors $\{ iso_V | V \in \mathcal{S} \}$ is a set of
projective generators of $\Gq$, where $\mathcal{S}$ is a set of
representatives of isometry classes of degenerate quadratic spaces. By
proposition \ref{4.32} $(2)$ $iso_V$ is a direct summand of 
$Q_V$; as $Q_V$ is of finite type, we deduce that $iso_V$ is of finite
type. Consequently, by the previous result we obtain that
$$\Gq \simeq \mathrm{Func}^{add}(\mathcal{P}^{op}, \mathcal{A}b)$$
where $\mathcal{P}$ is the full subcategory of $\Gq$ having as
objects the isotropic functors. By corollary \ref{4.33},
$\mathrm{Hom}_{\Gq} (iso_V, iso_W)=0$ if $V \not \simeq
W$. Consequently 
$$\mathrm{Func}^{add}(\mathcal{P}^{op}, \mathcal{A}b) \simeq \prod_{V \in
  O(V)} \mathrm{Func}^{add}(\mathcal{I}so_V^{op}, \mathcal{A}b)$$
where $\mathcal{I}so_V$ is the full subcategory of $\Gq$ with one
object, the functor $iso_V$. By corollary \ref{4.33}, $\mathrm{Hom}_{\Gq} (iso_V, iso_V)=\FF[O(V)]$; we deduce that  
$$\Gq \simeq  \prod_{V \in \mathcal{S}} \FF[O(V)]-mod.$$
\end{proof}
\appendix
\section{Properties of the precomposition functor}\label{appendix}
In this appendix we list several results about functor
categories. We have chosen to provide proofs since, even if these results are well-known, they are not easy to
find in the literature.

We are interested in the following question:

Let $\C$ and $\D$ be two categories, $\A$ be an abelian category,
$F: \C \rightarrow \D$ be a functor and $- \circ F: \mathrm{Func}(\D, \A)
\rightarrow \mathrm{Func}(\C, \A)$ be the precomposition functor, where
$\mathrm{Func}(\C, \A)$ is the category of functors from $\C$ to $\A$.

When the functor $F$ has a property $\mathcal{P}$, what can we deduce
for the precomposition functor?

Before giving three answers in the following propositions, we recall
that, by  \cite{gabriel}, $\mathrm{Func}(\C, \A)$ and $\mathrm{Func}(\D, \A)$
are abelian categories as the category $\A$ is abelian. Furthermore,
we remark that the precomposition functor is exact.
\begin{prop} \label{A1}
If $F$ is essentially surjective, then $- \circ F$ is faithful.
\end{prop}
\begin{proof}
As the precomposition functor is exact, it is sufficient to prove that, if
$H$ is an object of $\mathrm{Func}(\D, \A)$ such that $H \circ F=0$, then
$H=0$. For an object $D$ of $\mathcal{D}$, there exists an object $C$
of $\mathcal{C}$ such that $F(C) \simeq D$ as $F$ is essentially
surjective. So 
$$H(D) \simeq H(F(C)) =H \circ F(C)=0.$$

\end{proof}
\begin{prop} \label{A2}
If $F$ is full and essentially surjective, then:

\begin{enumerate}
\item 
the precomposition functor is fully-faithful;
\item
any subobject of an object in the image of the precomposition functor
is isomorphic to an object in the image of the precomposition functor;
\item
the image by the precomposition functor of a simple functor of
$\mathrm{Func}(\D, \A)$ is a simple functor of $\mathrm{Func}(\C, \A)$.
\end{enumerate}
\end{prop}
\begin{proof}
\begin{enumerate}

\item
According to proposition \ref{A1}, the functor $- \circ F$ is faithful.

For the fullness, we consider two objects $G$ and $H$ of
$\mathrm{Func}(\D, \A)$ and $\alpha$ an element of 
$\mathrm{Hom}_{\mathrm{Func}(\C, \A)}(G \circ F, H \circ F)$. We want
to prove that there exists a morphism $\beta$ of
$\mathrm{Hom}_{\mathrm{Func}(\D, \A)}(G, H)$ such that $\beta \circ F
=\alpha$.

Let $D$ be an object of $\D$, as $F$ is essentially surjective, we can
chose an object $C$ of $\C$ such that there exists an isomorphism:
$$\phi: F(C) \rightarrow D.$$
We define an element $\beta(D)$ of $\mathrm{Hom}_{\A}(G(D), H(D))$ by
the following composition 
$$G(D) \xrightarrow{G(\phi^{-1})} G \circ F(C) \xrightarrow{\alpha(C)}
H \circ F(C) \xrightarrow{H(\phi)} H(D).$$
Now, we show that $\beta$ defines a natural transformation from $G$ to $H$.

Let $f: D \rightarrow D'$ be an element of $\mathrm{Hom}_{\D}(D,D')$
and $\phi: F(C) \rightarrow D$ and $\phi': F(C') \rightarrow D'$ be the
isomorphisms associated to the choices of $C$ and $C'$ by the
essential surjectivity of $F$. As the functor $F$ is full, there
exists an element $g$ of $\mathrm{Hom}_{\C}(C,C')$ such that 
\begin{equation} \label{ann-eq1}
F(g)=\phi'^{-1} \circ f \circ \phi.
\end{equation}
Moreover, the following diagram is commutative
$$\xymatrix{
G(D) \ar[r]^{\beta(D)} \ar[d]_{G(\phi^{-1})}& H(D) \ar[d]^{H(\phi^{-1})} \\
G\circ F(C) \ar[r]^{\alpha(C)} \ar[d]_{G \circ F(g)}& H \circ F(C)
\ar[d]^{H\circ F(g)} \\
G\circ F(C') \ar[r]^{\alpha(C')} \ar[d]_{G(\phi')}& H \circ F(C')
\ar[d]^{H (\phi')} \\
G(D') \ar[r]^{\beta(D')} & H(D')  \\
}$$
because the higher square (respectively lower) commute by the
definition of $\beta(D)$ (respectively  $\beta(D')$) and the
commutativity of the square in
the center is a consequence of the naturality of
$\alpha$. Since 
$$G(\phi') \circ G F(g) \circ G(\phi^{-1})=G(\phi' \circ  F(g) \circ \phi^{-1})=G(f)$$
and 
$$H(\phi') \circ H F(g) \circ H(\phi^{-1})=H(\phi' \circ  F(g) \circ
\phi^{-1})=H(f)$$
where the first equalities rise from the functoriality of $G$ and
$H$ and the second ones rise from the relation \ref{ann-eq1}, we
deduce that $\beta$ is natural.
\item
Let $G$ be an object of $\mathrm{Func}(\D, \A)$ and $H$ be a subobject of
$G \circ F$. As in the first point, for an object $D$ of $\D$, by the
essential surjectivity of $F$, we can chose an object $C$ of $\C$ such
that there exists an isomorphism:
$$\phi: F(C) \rightarrow D.$$
We define an object of $\A$ by:
$$K(D):=H(C).$$

As in the first point, for $f: D \rightarrow D'$ a morphism of
$\mathrm{Hom}_{\D}(D,D')$, we denote by $\phi: F(C) \rightarrow D$ and $\phi': F(C') \rightarrow D'$ the
isomorphisms associated to the choices of $C$ and $C'$ by the
essential surjectivity of $F$. Since the functor $F$ is full, there
exists an element $g$ of $\mathrm{Hom}_{\C}(C,C')$ such that 
\begin{equation} \label{ann-eq2}
F(g)=\phi'^{-1} \circ f \circ \phi.
\end{equation}
We define a morphism
$$K(f): K(D) \rightarrow K(D')$$
by $K(f)=H(g)$. We obtain the following commutative diagram:
$$\xymatrix{
K(D) \ar@{=}[r] \ar[d]_{K(f)} & H(C) \ar@{^{(}->}[r] \ar[d]_{H(g)}& G
\circ F(C) \ar[r]^{G(\phi)} \ar[d]_{G \circ F(g)} &
G(D) \ar[d]_{G(f)}\\
K(D') \ar@{=}[r]  & H(C') \ar@{^{(}->}[r] & G
\circ F(C') \ar[r]^{G(\phi')} &
G(D').
}$$

Since the horizontal arrows of the diagram are monomorphisms, the
definition of $K(f)$ does not depend of the choice of the morphism $g$ in
\ref{ann-eq2}.
The functor $K$, thus defined, satisfies $K \circ F \simeq H$. Moreover,
we deduce from the commutativity of the above diagram that $K$ is
a subfunctor of $G$.

\item
Let $S$ be a simple functor of $\mathrm{Func}(\D, \A)$ and  $(- \circ F)(S)=S
\circ F$ its image by $- \circ F$. Suppose that $S \circ F$ is not
simple, so that there exists a proper subfunctor $H$ of $S \circ F$. By the proof
of the point $(2)$, of the proposition, there exists a subfunctor $K$
of $S$ such that:
$$H \simeq K \circ F.$$
Since $H$ is a proper subfunctor of $S \circ F$, $K$ is a proper
subfunctor of $S$. This contradicts the simplicity of $S$.

\end{enumerate}
\end{proof}
\begin{rem}
The constructions used in the proposition are closely related to the
left Kan extension. We have prefered to give an explicit proof rather
than a slicker one using the Kan extension.
\end{rem}

In the following proposition, we take $\mathcal{A}=Mod_A$ the category
of right $A$-modules where $A$ is a commutative ring with unit. By the
Yoneda lemma, we know that the functor
$$P_X^{\C}: \C \xrightarrow{\mathrm{Hom}_{\C}(X,-)} \mathcal{E}ns
\xrightarrow{A[-]} Mod_A$$
is a projective object of $\mathrm{Func}(\C, Mod_A)$, where $A[-]$ is the
linearisation functor. We have, in this case, the following
proposition.
\begin{prop} \label{A3}
If $F: \C \rightarrow \D$ is full and, for all objects $D$ of $\D$, there exists an
object $C$ of $\C$ such that $P_D^{\D}$ is a direct summand of $P_{FC}^{\D}$
then:
\begin{enumerate}
\item
the precomposition functor $- \circ F$ is faithful;
\item
a morphism $\sigma$ of $\mathrm{Func}(\D,Mod_A)$ such that $(- \circ
F)(\sigma)$ is an isomorphism, is an isomorphism;
\item
the categories $\mathrm{Func}(\D,Mod_A)$ and $\mathrm{Func}(\C,Mod_A)$ are equivalent.
\end{enumerate}
\end{prop}
\begin{proof}
\begin{enumerate}
\item
As the precomposition functor is exact, it is sufficient to prove that, for
$H$ an object of $\mathrm{Func}(\D, \A)$ such that $H \circ F=0$,
$H=0$. Let $D$ be an object of $\D$, then $H(D)=Hom(P_D^{\D},H)$ by the Yoneda
lemma. By hypothesis, this is a direct summand of 
$$Hom(P_{FC}^{\D},H) \simeq H \circ F(C)=0.$$
Consequently $H(D)=0$.

\item
Let $\sigma: F \rightarrow G$ be a morphism of $\mathrm{Func}(\D, Mod_A)$. We
have the following exact sequence:
$$0 \rightarrow Ker(\sigma) \rightarrow F \xrightarrow{\sigma} G
\rightarrow Coker(\sigma) \rightarrow 0.$$
As $(- \circ H)(\sigma)$ is an isomorphism, $Ker(\sigma)
\circ H=0$ and $Coker(\sigma) \circ H =0$. Consequently,
$Ker(\sigma)=Coker(\sigma)=0$ by the argument of the previous point. Hence,
$\sigma$ is an isomorphism.

\item
By general results, the functor $(- \circ F): \mathrm{Func}(\D,Mod_A)
\rightarrow \mathrm{Func}(\C,Mod_A)$  admits a right adjoint, given by Kan
extension. We will denote by $(\tilde{-})$ this adjoint. We will prove
that $(\tilde{-})$ and $(- \circ F)$ define an equivalence of
categories. 

As $F$ is full, we have, for all objects $W$ of $\C$
$$P_{FC}^{\D}(FW)=A[Hom_{\D}(FC,FW)]=A[Hom_{\C}(C,W)]=P_C^{\C}(W).$$
It follows that $P_{FC}^{\D} \circ F=P_C^{\C}$. As $P_D^{\D} \circ F$ is a direct summand of
  $P_{FC}^{\D} \circ F$ by hypothesis, we deduce that $P_D^{\D} \circ F$ is a
    direct summand of $P_C^{\D}$ and hence projective.
Consequently, $(\tilde{-})$ is an exact functor since
$$\tilde{H}(D) \simeq \mathrm{Hom}(P_D^{\D}, \tilde{H}) \simeq
\mathrm{Hom}(P_D^{\D} \circ F, H).$$

Let $H$ be an object of $\mathrm{Func}(\C,Mod_A)$, then there is a
sequence of isomorphisms:
$$(- \circ F) \circ (\tilde{-})(H)=\tilde{H}(FZ) \simeq
Hom(P_{FZ}^{\D} \circ F, H) \simeq Hom(P_Z^{\C},H) \simeq H(Z).$$
So, the counit of the adjunction
\begin{equation} \label{eq-ann-3}
(- \circ F) \circ (\tilde{-}) \rightarrow Id
\end{equation}
is an isomorphism.

Consequently, the unit of the adjunction
\begin{equation} 
Id  \rightarrow  (\tilde{-}) \circ (- \circ F) 
\end{equation}
induces, by composition with $(- \circ F)$, an isomorphism
\begin{equation} 
(- \circ F)  \rightarrow (- \circ F) \circ (\tilde{-}) \circ (- \circ F). 
\end{equation}

So, by the previous point of the proposition, the unit  of the
adjunction $ Id \rightarrow (\tilde{-})
\circ (- \circ F)$ is an isomorphism.

\end{enumerate}
\end{proof}


\bibliographystyle{amsplain}
\bibliography{these}

\end{document}